\documentclass{amsart}
\usepackage{amsfonts}

\newtheorem{theorem}{Theorem}
\theoremstyle{plain}

\newtheorem{corollary}{Corollary}

\newtheorem{proposition}{Proposition}
\newtheorem{remark}{Remark}

\numberwithin{equation}{section}
\input{tcilatex}

\begin{document}
\title[Reverses of Schwarz, Triangle and Bessel Inequalities]{Reverses of
Schwarz, Triangle and Bessel Inequalities in Inner Product Spaces}
\author{S.S. Dragomir}
\address{School of Computer Science and Mathematics\\
Victoria University of Technology\\
PO Box 14428, MCMC 8001\\
Victoria, Australia.}
\email{sever.dragomir@vu.edu.au}
\urladdr{http://rgmia.vu.edu.au/SSDragomirWeb.html}
\date{August 04, 2003.}
\subjclass[2000]{Primary 26D15, 46C05.}
\keywords{Schwarz's inequality, Triangle inequality, Bessel's inequality, Gr%
\"{u}ss type inequalities, Integral inequalities.}

\begin{abstract}
Reverses of Schwarz, triangle and Bessel inequalities in inner product
spaces that improve some earlier results are pointed out. They are applied
to obtain new Gr\"{u}ss type inequalities in inner product spaces. Some
natural applications for integral inequalities are also pointed out.
\end{abstract}

\maketitle

\section{Introduction\label{s1}}

Let $\left( H;\left\langle \cdot ,\cdot \right\rangle \right) $ be an inner
product over the real or complex number field $\mathbb{K}$. The following
inequality is known in the literature as \textit{Schwarz's inequality}:%
\begin{equation}
\left\vert \left\langle x,y\right\rangle \right\vert ^{2}\leq \left\Vert
x\right\Vert ^{2}\left\Vert y\right\Vert ^{2},\ \ \ \ x,y\in H;  \label{1.1}
\end{equation}%
where $\left\Vert z\right\Vert ^{2}=\left\langle z,z\right\rangle ,$ $z\in
H. $ The equality occurs in (\ref{1.1}) if and only if $x$ and $y$ are
linearly dependent.

In \cite{SSD1}, the following \textit{reverse} of Schwarz's inequality has
been obtained:%
\begin{equation}
0\leq \left\Vert x\right\Vert ^{2}\left\Vert y\right\Vert ^{2}-\left\vert
\left\langle x,y\right\rangle \right\vert ^{2}\leq \frac{1}{4}\left\vert
A-a\right\vert ^{2}\left\Vert y\right\Vert ^{4},  \label{1.2}
\end{equation}%
provided $x,y\in H$ and $a,A\in \mathbb{K}$ are so that either%
\begin{equation}
\func{Re}\left\langle Ay-x,x-ay\right\rangle \geq 0,  \label{1.3}
\end{equation}%
or, equivalently,%
\begin{equation}
\left\Vert x-\frac{a+A}{2}\cdot y\right\Vert \leq \frac{1}{2}\left\vert
A-a\right\vert \left\Vert y\right\Vert ,  \label{1.4}
\end{equation}%
holds. The constant $\frac{1}{4}$ is best possible in (\ref{1.2}) in the
sense that it cannot be replaced by a smaller constant.

If $x,y,A,a$ satisfy either (\ref{1.3}) or (\ref{1.4}), then the following
reverse of Schwarz's inequality also holds \cite{SSD2}%
\begin{align}
\left\Vert x\right\Vert \left\Vert y\right\Vert & \leq \frac{1}{2}\cdot 
\frac{\func{Re}\left[ A\overline{\left\langle x,y\right\rangle }+\overline{a}%
\left\langle x,y\right\rangle \right] }{\left[ \func{Re}\left( \overline{a}%
A\right) \right] ^{\frac{1}{2}}}  \label{1.5} \\
& \leq \frac{1}{2}\cdot \frac{\left\vert A\right\vert +\left\vert
a\right\vert }{\left[ \func{Re}\left( \overline{a}A\right) \right] ^{\frac{1%
}{2}}}\left\vert \left\langle x,y\right\rangle \right\vert ,  \notag
\end{align}%
provided that, the complex numbers $a$ and $A$ satisfy the condition $\func{%
Re}\left( \overline{a}A\right) >0.$ In both inequalities in (\ref{1.5}), the
constant $\frac{1}{2}$ is best possible.

An additive version of (\ref{1.5}) may be stated as well (see also \cite%
{SSD3})%
\begin{equation}
0\leq \left\Vert x\right\Vert ^{2}\left\Vert y\right\Vert ^{2}-\left\vert
\left\langle x,y\right\rangle \right\vert ^{2}\leq \frac{1}{4}\cdot \frac{%
\left( \left\vert A\right\vert -\left\vert a\right\vert \right) ^{2}+4\left[
\left\vert Aa\right\vert -\func{Re}\left( \overline{a}A\right) \right] }{%
\func{Re}\left( \overline{a}A\right) }\left\vert \left\langle
x,y\right\rangle \right\vert ^{2}.  \label{1.6}
\end{equation}%
In this inequality, $\frac{1}{4}$ is the best possible constant.

It has been proven in \cite{SSD4}, that%
\begin{equation}
0\leq \left\Vert x\right\Vert ^{2}-\left\vert \left\langle x,y\right\rangle
\right\vert ^{2}\leq \frac{1}{4}\left\vert \phi -\varphi \right\vert
^{2}-\left\vert \frac{\phi +\varphi }{2}-\left\langle x,e\right\rangle
\right\vert ^{2};  \label{1.7}
\end{equation}%
provided, either 
\begin{equation}
\func{Re}\left\langle \phi e-x,x-\varphi e\right\rangle \geq 0,  \label{1.8}
\end{equation}%
or, equivalently,%
\begin{equation}
\left\Vert x-\frac{\phi +\varphi }{2}e\right\Vert \leq \frac{1}{2}\left\vert
\phi -\varphi \right\vert ,  \label{1.9}
\end{equation}%
where $e=H,$ $\left\Vert e\right\Vert =1.$ The constant $\frac{1}{4}$ in \ref%
{1.7} is also best possible.

If we choose $e=\frac{y}{\left\Vert y\right\Vert },$ $\phi =\Gamma
\left\Vert y\right\Vert ,$ $\varphi =\gamma \left\Vert y\right\Vert $ $%
\left( y\neq 0\right) ,$ $\Gamma ,\gamma \in \mathbb{K}$, then by (\ref{1.8}%
), (\ref{1.9}) we have,%
\begin{equation}
\func{Re}\left\langle \Gamma y-x,x-\gamma y\right\rangle \geq 0,
\label{1.10}
\end{equation}%
or, equivalently,%
\begin{equation}
\left\Vert x-\frac{\Gamma +\gamma }{2}y\right\Vert \leq \frac{1}{2}%
\left\vert \Gamma -\gamma \right\vert \left\Vert y\right\Vert ,  \label{1.11}
\end{equation}%
imply the following reverse of Schwarz's inequality:%
\begin{equation}
0\leq \left\Vert x\right\Vert ^{2}\left\Vert y\right\Vert ^{2}-\left\vert
\left\langle x,y\right\rangle \right\vert ^{2}\leq \frac{1}{4}\left\vert
\Gamma -\gamma \right\vert ^{2}\left\Vert y\right\Vert ^{4}-\left\vert \frac{%
\Gamma +\gamma }{2}\left\Vert y\right\Vert ^{2}-\left\langle
x,y\right\rangle \right\vert ^{2}.  \label{1.12}
\end{equation}%
The constant $\frac{1}{4}$ in (\ref{1.12}) is sharp.

Note that this inequality is an improvement of (\ref{1.2}), but it might not
be very convenient for applications.

Now, let $\left\{ e_{i}\right\} _{i\in I}$ be a finite or infinite family of
orthornormal vectors in the inner product space $\left( H;\left\langle \cdot
,\cdot \right\rangle \right) ,$ i.e., we recall that 
\begin{equation*}
\left\langle e_{i},e_{j}\right\rangle =\left\{ 
\begin{array}{ll}
0 & \text{if \ }i\neq j \\ 
&  \\ 
1 & \text{if \ }i=j%
\end{array}%
\right. ,\ \ \ i,j\in I.
\end{equation*}%
In \cite{SSD5}, we proved that if $\left\{ e_{i}\right\} _{i\in I}$ is as
above, $F\subset I$ is a finite part of $I$ such that either%
\begin{equation}
\func{Re}\left\langle \sum_{i\in F}\phi _{i}e_{i}-x,x-\sum_{i\in F}\varphi
_{i}e_{i}\right\rangle \geq 0,  \label{1.13}
\end{equation}%
or, equivalently,%
\begin{equation}
\left\Vert x-\sum_{i\in F}\frac{\phi _{i}+\varphi _{i}}{2}e_{i}\right\Vert
\leq \frac{1}{2}\left( \sum_{i\in F}\left\vert \phi _{i}-\varphi
_{i}\right\vert ^{2}\right) ^{\frac{1}{2}},  \label{1.14}
\end{equation}%
holds, where $\left( \phi _{i}\right) _{i\in I},$ $\left( \varphi
_{i}\right) _{i\in I}$ are real or complex numbers, then we have the
following reverse of \textit{Bessel's inequality:}%
\begin{align}
0& \leq \left\Vert x\right\Vert ^{2}-\sum_{i\in F}\left\vert \left\langle
x,e_{i}\right\rangle \right\vert ^{2}  \label{1.15} \\
& \leq \frac{1}{4}\cdot \sum_{i\in F}\left\vert \phi _{i}-\varphi
_{i}\right\vert ^{2}-\func{Re}\left\langle \sum_{i\in F}\phi
_{i}e_{i}-x,x-\sum_{i\in F}\varphi _{i}e_{i}\right\rangle  \notag \\
& \leq \frac{1}{4}\cdot \sum_{i\in F}\left\vert \phi _{i}-\varphi
_{i}\right\vert ^{2}.  \notag
\end{align}%
The constant $\frac{1}{4}$ in both inequalities is sharp. This result
improves an earlier result by N. Ujevi\'{c} obtained only for real spaces 
\cite{NU}.

In \cite{SSD4}, by the use of a different technique, another reverse of
Bessel's inequality has been proven, namely:%
\begin{align}
0& \leq \left\Vert x\right\Vert ^{2}-\sum_{i\in F}\left\vert \left\langle
x,e_{i}\right\rangle \right\vert ^{2}  \label{1.16} \\
& \leq \frac{1}{4}\cdot \sum_{i\in F}\left\vert \phi _{i}-\varphi
_{i}\right\vert ^{2}-\sum_{i\in F}\left\vert \frac{\phi _{i}+\varphi _{i}}{2}%
-\left\langle x,e_{i}\right\rangle \right\vert ^{2}  \notag \\
& \leq \frac{1}{4}\cdot \sum_{i\in F}\left\vert \phi _{i}-\varphi
_{i}\right\vert ^{2},  \notag
\end{align}%
provided that $\left( e_{i}\right) _{i\in I},$ $\left( \phi _{i}\right)
_{i\in I},$ $\left( \varphi _{i}\right) _{i\in I},$ $x$ and $F$ are as above.

Here the constant $\frac{1}{4}$ is sharp in both inequalities.

It has also been shown that the bounds provided by (\ref{1.15}) and (\ref%
{1.16}) for the Bessel's difference $\left\Vert x\right\Vert ^{2}-\sum_{i\in
F}\left\vert \left\langle x,e_{i}\right\rangle \right\vert ^{2}$ cannot be
compared in general, meaning that there are examples for which one is
smaller than the other \cite{SSD4}.

Finally, we recall another type of reverse for Bessel inequality that has
been obtained in \cite{SSD6}:%
\begin{equation}
\left\Vert x\right\Vert ^{2}\leq \frac{1}{4}\cdot \frac{\sum_{i\in F}\left(
\left\vert \phi _{i}\right\vert +\left\vert \varphi _{i}\right\vert \right)
^{2}}{\sum_{i\in F}\func{Re}\left( \phi _{i}\overline{\varphi _{i}}\right) }%
\sum_{i\in F}\left\vert \left\langle x,e_{i}\right\rangle \right\vert ^{2};
\label{1.17}
\end{equation}%
provided $\left( \phi _{i}\right) _{i\in I},$ $\left( \varphi _{i}\right)
_{i\in I}$ satisfy (\ref{1.13}) (or, equivalently (\ref{1.14})) and $%
\sum_{i\in F}\func{Re}\left( \phi _{i}\overline{\varphi _{i}}\right) >0.$
Here the constant $\frac{1}{4}$ is also best possible.

An additive version of (\ref{1.17}) is 
\begin{align}
0& \leq \left\Vert x\right\Vert ^{2}-\sum_{i\in F}\left\vert \left\langle
x,e_{i}\right\rangle \right\vert ^{2}  \label{1.18} \\
& \leq \frac{1}{4}\cdot \frac{\sum_{i\in F}\left\{ \left( \left\vert \phi
_{i}\right\vert -\left\vert \varphi _{i}\right\vert \right) ^{2}+4\left[
\left\vert \phi _{i}\varphi _{i}\right\vert -\func{Re}\left( \phi _{i}%
\overline{\varphi _{i}}\right) \right] \right\} }{\sum_{i\in F}\func{Re}%
\left( \phi _{i}\overline{\varphi _{i}}\right) }.  \notag
\end{align}%
The constant $\frac{1}{4}$ is best possible.

It is the main aim of the present paper to point out new reverse
inequalities to Schwarz's, triangle and Bessel's inequalities.

Some results related to Gr\"{u}ss' inequality in inner product spaces are
also pointed out. Natural applications for integrals are also provided.

\section{Some Reverses of Schwarz's Inequality\label{s2}}

The following result holds.

\begin{theorem}
\label{t2.1}Let $\left( H;\left\langle \cdot ,\cdot \right\rangle \right) $
be an inner product space over the real or complex number field $\mathbb{K}$ 
$\left( \mathbb{K}=\mathbb{R},\ \mathbb{K}=\mathbb{C}\right) $ and $x,a\in
H, $ $r>0$ are such that%
\begin{equation}
x\in \overline{B}\left( a,r\right) :=\left\{ z\in H|\left\Vert
z-a\right\Vert \leq r\right\} .  \label{2.1}
\end{equation}

\begin{enumerate}
\item[(i)] If $\left\Vert a\right\Vert >r,$ then we have the inequality%
\begin{equation}
0\leq \left\Vert x\right\Vert ^{2}\left\Vert a\right\Vert ^{2}-\left\vert
\left\langle x,a\right\rangle \right\vert ^{2}\leq \left\Vert x\right\Vert
^{2}\left\Vert a\right\Vert ^{2}-\left[ \func{Re}\left\langle
x,a\right\rangle \right] ^{2}\leq r^{2}\left\Vert x\right\Vert ^{2}.
\label{2.2}
\end{equation}%
The constant $C=1$ in front of $r^{2}$ is best possible in the sense that it
cannot be replaced by a smaller one.

\item[(ii)] If $\left\Vert a\right\Vert =r,$ then%
\begin{equation}
\left\Vert x\right\Vert ^{2}\leq 2\func{Re}\left\langle x,a\right\rangle
\leq 2\left\vert \left\langle x,a\right\rangle \right\vert .  \label{2.3}
\end{equation}%
The constant $2$ is best possible in both inequalities.

\item[(iii)] If $\left\Vert a\right\Vert <r,$ then%
\begin{equation}
\left\Vert x\right\Vert ^{2}\leq r^{2}-\left\Vert a\right\Vert ^{2}+2\func{Re%
}\left\langle x,a\right\rangle \leq r^{2}-\left\Vert a\right\Vert
^{2}+2\left\vert \left\langle x,a\right\rangle \right\vert .  \label{2.4}
\end{equation}%
Here the constant $2$ is also best possible.
\end{enumerate}
\end{theorem}

\begin{proof}
Since $x\in \overline{B}\left( a,r\right) ,$ then obviously $\left\Vert
x-a\right\Vert ^{2}\leq r^{2},$ which is equivalent to 
\begin{equation}
\left\Vert x\right\Vert ^{2}+\left\Vert a\right\Vert ^{2}-r^{2}\leq 2\func{Re%
}\left\langle x,a\right\rangle .  \label{2.5}
\end{equation}

\begin{enumerate}
\item[(i)] If $\left\Vert a\right\Vert >r,$ then we may divide (\ref{2.5})
by $\sqrt{\left\Vert a\right\Vert ^{2}-r^{2}}>0$ getting 
\begin{equation}
\frac{\left\Vert x\right\Vert ^{2}}{\sqrt{\left\Vert a\right\Vert ^{2}-r^{2}}%
}+\sqrt{\left\Vert a\right\Vert ^{2}-r^{2}}\leq \frac{2\func{Re}\left\langle
x,a\right\rangle }{\sqrt{\left\Vert a\right\Vert ^{2}-r^{2}}}.  \label{2.6}
\end{equation}%
Using the elementary inequality%
\begin{equation*}
\alpha p+\frac{1}{\alpha }q\geq 2\sqrt{pq},\ \ \ \alpha >0,\ \ p,q\geq 0,
\end{equation*}%
we may state that%
\begin{equation}
2\left\Vert x\right\Vert \leq \frac{\left\Vert x\right\Vert ^{2}}{\sqrt{%
\left\Vert a\right\Vert ^{2}-r^{2}}}+\sqrt{\left\Vert a\right\Vert ^{2}-r^{2}%
}.  \label{2.7}
\end{equation}%
Making use of (\ref{2.6}) and (\ref{2.7}), we deduce%
\begin{equation}
\left\Vert x\right\Vert \sqrt{\left\Vert a\right\Vert ^{2}-r^{2}}\leq \func{%
Re}\left\langle x,a\right\rangle .  \label{2.8}
\end{equation}%
Taking the square in (\ref{2.8}) and re-arranging the terms, we deduce the
third inequality in (\ref{2.2}). The others are obvious.

To prove the sharpness of the constant, assume, under the hypothesis of the
theorem, that, there exists a constant $c>0$ such that%
\begin{equation}
\left\Vert x\right\Vert ^{2}\left\Vert a\right\Vert ^{2}-\left[ \func{Re}%
\left\langle x,a\right\rangle \right] ^{2}\leq cr^{2}\left\Vert x\right\Vert
^{2},  \label{2.9}
\end{equation}%
provided $x\in \overline{B}\left( a,r\right) $ and $\left\Vert a\right\Vert
>r.$

Let $r=\sqrt{\varepsilon }>0,$ $\varepsilon \in \left( 0,1\right) ,$ $a,e\in
H$ with $\left\Vert a\right\Vert =\left\Vert e\right\Vert =1$ and $a\perp e.$
Put $x=a+\sqrt{\varepsilon }e.$ Then obviously $x\in \overline{B}\left(
a,r\right) ,$ $\left\Vert a\right\Vert >r$ and $\left\Vert x\right\Vert
^{2}=\left\Vert a\right\Vert ^{2}+\varepsilon \left\Vert e\right\Vert
^{2}=1+\varepsilon $, $\func{Re}\left\langle x,a\right\rangle =\left\Vert
a\right\Vert ^{2}=1,$ and thus $\left\Vert x\right\Vert ^{2}\left\Vert
a\right\Vert ^{2}-\left[ \func{Re}\left\langle x,a\right\rangle \right]
^{2}=\varepsilon .$ Using (\ref{2.9}), we may write that%
\begin{equation*}
\varepsilon \leq c\varepsilon \left( 1+\varepsilon \right) ,\ \ \varepsilon
>0
\end{equation*}%
giving 
\begin{equation}
c+c\varepsilon \geq 1\text{ \ for any }\varepsilon >0  \label{2.10}
\end{equation}%
Letting $\varepsilon \rightarrow 0+,$ we get from (\ref{2.10}) that $c\geq 1,
$ and the sharpness of the constant is proved.

\item[(ii)] The inequality (\ref{2.3}) is obvious by (\ref{2.5}) since $%
\left\Vert a\right\Vert =r.$ The best constant follows in a similar way to
the above.

\item[(iii)] The inequality (\ref{2.3}) is obvious. The best constant may be
proved in a similar way to the above. We omit the details.
\end{enumerate}
\end{proof}

The following reverse of Schwarz's inequality holds.

\begin{theorem}
\label{t2.2}Let $\left( H;\left\langle \cdot ,\cdot \right\rangle \right) $
be an inner product space over $\mathbb{K}$ and $x,y\in H,$ $\gamma ,\Gamma
\in \mathbb{K}$ such that either%
\begin{equation}
\func{Re}\left\langle \Gamma y-x,x-\gamma y\right\rangle \geq 0,
\label{2.11}
\end{equation}%
or, equivalently,%
\begin{equation}
\left\Vert x-\frac{\Gamma +\gamma }{2}y\right\Vert \leq \frac{1}{2}%
\left\vert \Gamma -\gamma \right\vert \left\Vert y\right\Vert ,  \label{2.12}
\end{equation}%
holds.

\begin{enumerate}
\item[(i)] If $\func{Re}\left( \Gamma \overline{\gamma }\right) >0,$ then we
have the inequalities%
\begin{align}
\left\Vert x\right\Vert ^{2}\left\Vert y\right\Vert ^{2}& \leq \frac{1}{4}%
\cdot \frac{\left\{ \func{Re}\left[ \left( \overline{\Gamma }+\overline{%
\gamma }\right) \left\langle x,y\right\rangle \right] \right\} ^{2}}{\func{Re%
}\left( \Gamma \overline{\gamma }\right) }  \label{2.13} \\
& \leq \frac{1}{4}\cdot \frac{\left\vert \Gamma +\gamma \right\vert ^{2}}{%
\func{Re}\left( \Gamma \overline{\gamma }\right) }\left\vert \left\langle
x,y\right\rangle \right\vert ^{2}.  \notag
\end{align}%
The constant $\frac{1}{4}$ is best possible in both inequalities.

\item[(ii)] If $\func{Re}\left( \Gamma \overline{\gamma }\right) =0,$ then 
\begin{equation}
\left\Vert x\right\Vert ^{2}\leq \func{Re}\left[ \left( \overline{\Gamma }+%
\overline{\gamma }\right) \left\langle x,y\right\rangle \right] \leq
\left\vert \Gamma +\gamma \right\vert \left\vert \left\langle
x,y\right\rangle \right\vert .  \label{2.14}
\end{equation}

\item[(iii)] If $\func{Re}\left( \Gamma \overline{\gamma }\right) <0,$ then 
\begin{align}
\left\Vert x\right\Vert ^{2}& \leq -\func{Re}\left( \Gamma \overline{\gamma }%
\right) \left\Vert y\right\Vert ^{2}+\func{Re}\left[ \left( \overline{\Gamma 
}+\overline{\gamma }\right) \left\langle x,y\right\rangle \right]
\label{2.15} \\
& \leq -\func{Re}\left( \Gamma \overline{\gamma }\right) \left\Vert
y\right\Vert ^{2}+\left\vert \Gamma +\gamma \right\vert \left\vert
\left\langle x,y\right\rangle \right\vert .  \notag
\end{align}
\end{enumerate}
\end{theorem}

\begin{proof}
The proof of the equivalence between the inequalities (\ref{2.11}) and (\ref%
{2.12}) follows by the fact that in an inner product space $\func{Re}%
\left\langle Z-x,x-z\right\rangle \geq 0$ for $x,z,Z\in H$ is equivalent
with $\left\Vert x-\frac{z+Z}{2}\right\Vert \leq \frac{1}{2}\left\Vert
Z-z\right\Vert $ (see for example \cite{SSD3}).

Consider, for $y\neq 0,$ $a=\frac{\gamma +\Gamma }{2}y$ and $r=\frac{1}{2}%
\left\vert \Gamma -\gamma \right\vert \left\Vert y\right\Vert ^{2}.$ Then%
\begin{equation*}
\left\Vert a\right\Vert ^{2}-r^{2}=\frac{\left\vert \Gamma +\gamma
\right\vert ^{2}-\left\vert \Gamma -\gamma \right\vert ^{2}}{4}\left\Vert
y\right\Vert ^{2}=\func{Re}\left( \Gamma \overline{\gamma }\right)
\left\Vert y\right\Vert ^{2}.
\end{equation*}

\begin{enumerate}
\item[(i)] If $\func{Re}\left( \Gamma \overline{\gamma }\right) >0,$ then
the hypothesis of (i) in Theorem \ref{t2.1} is satisfied, and by the second
inequality in (\ref{2.2}) we have%
\begin{equation*}
\left\Vert x\right\Vert ^{2}\frac{\left\vert \Gamma +\gamma \right\vert ^{2}%
}{4}\left\Vert y\right\Vert ^{2}-\frac{1}{4}\left\{ \func{Re}\left[ \left( 
\overline{\Gamma }+\overline{\gamma }\right) \left\langle x,y\right\rangle %
\right] \right\} ^{2}\leq \frac{1}{4}\left\vert \Gamma -\gamma \right\vert
^{2}\left\Vert x\right\Vert ^{2}\left\Vert y\right\Vert ^{2}
\end{equation*}%
from where we derive%
\begin{equation*}
\frac{\left\vert \Gamma +\gamma \right\vert ^{2}-\left\vert \Gamma -\gamma
\right\vert ^{2}}{4}\left\Vert x\right\Vert ^{2}\left\Vert y\right\Vert
^{2}\leq \frac{1}{4}\left\{ \func{Re}\left[ \left( \overline{\Gamma }+%
\overline{\gamma }\right) \left\langle x,y\right\rangle \right] \right\}
^{2},
\end{equation*}%
giving the first inequality in (\ref{2.13}).

The second inequality is obvious.

To prove the sharpness of the constant $\frac{1}{4},$ assume that the first
inequality in (\ref{2.13}) holds with a constant $c>0,$ i.e., 
\begin{equation}
\left\Vert x\right\Vert ^{2}\left\Vert y\right\Vert ^{2}\leq c\cdot \frac{%
\left\{ \func{Re}\left[ \left( \overline{\Gamma }+\overline{\gamma }\right)
\left\langle x,y\right\rangle \right] \right\} ^{2}}{\func{Re}\left( \Gamma 
\overline{\gamma }\right) },  \label{2.16}
\end{equation}%
provided $\func{Re}\left( \Gamma \overline{\gamma }\right) >0$ and either (%
\ref{2.11}) or (\ref{2.12}) holds.

Assume that $\Gamma ,\gamma >0,$ and let $x=\gamma y.$ Then (\ref{2.11})
holds and by (\ref{2.16}) we deduce%
\begin{equation*}
\gamma ^{2}\left\Vert y\right\Vert ^{4}\leq c\cdot \frac{\left( \Gamma
+\gamma \right) ^{2}\gamma ^{2}\left\Vert y\right\Vert ^{4}}{\Gamma \gamma }
\end{equation*}%
giving%
\begin{equation}
\Gamma \gamma \leq c\left( \Gamma +\gamma \right) ^{2}\text{ \ for any \ }%
\Gamma ,\gamma >0.  \label{2.17}
\end{equation}%
Let $\varepsilon \in \left( 0,1\right) $ and choose in (\ref{2.17}), $\Gamma
=1+\varepsilon ,$ $\gamma =1-\varepsilon >0$ to get $1-\varepsilon ^{2}\leq
4c$ for any $\varepsilon \in \left( 0,1\right) .$ Letting $\varepsilon
\rightarrow 0+,$ we deduce $c\geq \frac{1}{4},$ and the sharpness of the
constant is proved.

(ii) and (iii) are obvious and we omit the details.
\end{enumerate}
\end{proof}

\begin{remark}
We observe that the second bound in (\ref{2.13}) for $\left\Vert
x\right\Vert ^{2}\left\Vert y\right\Vert ^{2}$ is better than the second
bound provided by (\ref{1.5}).
\end{remark}

The following corollary provides a reverse inequality for the additive
version of Schwarz's inequality.

\begin{corollary}
\label{c2.3}With the assumptions of Theorem \ref{t2.2} and if $\func{Re}%
\left( \Gamma \overline{\gamma }\right) >0,$ then we have the inequality:%
\begin{equation}
0\leq \left\Vert x\right\Vert ^{2}\left\Vert y\right\Vert ^{2}-\left\vert
\left\langle x,y\right\rangle \right\vert ^{2}\leq \frac{1}{4}\cdot \frac{%
\left\vert \Gamma -\gamma \right\vert ^{2}}{\func{Re}\left( \Gamma \overline{%
\gamma }\right) }\left\vert \left\langle x,y\right\rangle \right\vert ^{2}.
\label{2.18}
\end{equation}%
The constant $\frac{1}{4}$ is best possible in (\ref{2.18}).
\end{corollary}

The proof is obvious from (\ref{2.13}) on subtracting in both sides the same
quantity $\left\vert \left\langle x,y\right\rangle \right\vert ^{2}.$ The
sharpness of the constant may be proven in a similar manner to the one
incorporated in the proof of (i), Theorem \ref{t2.2}. We omit the details.

\begin{remark}
It is obvious that the inequality (\ref{2.18}) is better than (\ref{1.6})
obtained in \cite{SSD3}.
\end{remark}

For some recent results in connection to Schwarz's inequality, see \cite{ADR}%
, \cite{DM} and \cite{GH}.

\section{Reverses of the Triangle Inequality\label{s3}}

The following reverse of the triangle inequality holds.

\begin{proposition}
\label{p2.4}Let $\left( H;\left\langle \cdot ,\cdot \right\rangle \right) $
be an inner product space over the real or complex number field $\mathbb{K}$ 
$\left( \mathbb{K}=\mathbb{R},\mathbb{C}\right) $ and $x,a\in H,$ $r>0$ are
such that 
\begin{equation}
\left\Vert x-a\right\Vert \leq r<\left\Vert a\right\Vert .  \label{2.19}
\end{equation}%
Then we have the inequality%
\begin{equation}
0\leq \left\Vert x\right\Vert +\left\Vert a\right\Vert -\left\Vert
x+a\right\Vert \leq \sqrt{2}r\cdot \sqrt{\frac{\func{Re}\left\langle
x,a\right\rangle }{\sqrt{\left\Vert a\right\Vert ^{2}-r^{2}}\left( \sqrt{%
\left\Vert a\right\Vert ^{2}-r^{2}}+\left\Vert a\right\Vert \right) }}.
\label{2.20}
\end{equation}
\end{proposition}

\begin{proof}
Using the inequality (\ref{2.8}), we may write that%
\begin{equation*}
\left\Vert x\right\Vert \left\Vert a\right\Vert \leq \frac{\left\Vert
a\right\Vert \func{Re}\left\langle x,a\right\rangle }{\sqrt{\left\Vert
a\right\Vert ^{2}-r^{2}}},
\end{equation*}%
giving%
\begin{align}
0& \leq \left\Vert x\right\Vert \left\Vert a\right\Vert -\func{Re}%
\left\langle x,a\right\rangle  \label{2.21} \\
& \leq \func{Re}\left\langle x,a\right\rangle \frac{\left\Vert a\right\Vert -%
\sqrt{\left\Vert a\right\Vert ^{2}-r^{2}}}{\sqrt{\left\Vert a\right\Vert
^{2}-r^{2}}}  \notag \\
& =\frac{r^{2}\func{Re}\left\langle x,a\right\rangle }{\sqrt{\left\Vert
a\right\Vert ^{2}-r^{2}}\left( \sqrt{\left\Vert a\right\Vert ^{2}-r^{2}}%
+\left\Vert a\right\Vert \right) }.  \notag
\end{align}%
Since%
\begin{equation*}
\left( \left\Vert x\right\Vert +\left\Vert a\right\Vert \right)
^{2}-\left\Vert x+a\right\Vert ^{2}=2\left( \left\Vert x\right\Vert
\left\Vert a\right\Vert -\func{Re}\left\langle x,a\right\rangle \right) ,
\end{equation*}%
then by (\ref{2.21}), we have%
\begin{align*}
\left\Vert x\right\Vert +\left\Vert a\right\Vert & \leq \sqrt{\left\Vert
x+a\right\Vert ^{2}+\frac{2r^{2}\func{Re}\left\langle x,a\right\rangle }{%
\sqrt{\left\Vert a\right\Vert ^{2}-r^{2}}\left( \sqrt{\left\Vert
a\right\Vert ^{2}-r^{2}}+\left\Vert a\right\Vert \right) }} \\
& \leq \left\Vert x+a\right\Vert +\sqrt{2}r\cdot \sqrt{\frac{\func{Re}%
\left\langle x,a\right\rangle }{\sqrt{\left\Vert a\right\Vert ^{2}-r^{2}}%
\left( \sqrt{\left\Vert a\right\Vert ^{2}-r^{2}}+\left\Vert a\right\Vert
\right) }},
\end{align*}%
giving the desired inequality (\ref{2.20}).
\end{proof}

The following proposition providing a simpler reverse for the triangle
inequality also holds.

\begin{proposition}
\label{p2.5}Let $\left( H;\left\langle \cdot ,\cdot \right\rangle \right) $
be an inner product space over $\mathbb{K}$ and $x,y\in H,$ $M>m>0$ such
that either%
\begin{equation}
\func{Re}\left\langle My-x,x-my\right\rangle \geq 0,  \label{2.22}
\end{equation}%
or, equivalently,%
\begin{equation}
\left\Vert x-\frac{M+m}{2}\cdot y\right\Vert \leq \frac{1}{2}\left(
M-m\right) \left\Vert y\right\Vert ,  \label{2.23}
\end{equation}%
holds. Then we have the inequality%
\begin{equation}
0\leq \left\Vert x\right\Vert +\left\Vert y\right\Vert -\left\Vert
x+y\right\Vert \leq \frac{\sqrt{M}-\sqrt{m}}{\sqrt[4]{mM}}\sqrt{\func{Re}%
\left\langle x,y\right\rangle }.  \label{2.24}
\end{equation}
\end{proposition}

\begin{proof}
Choosing in (\ref{2.8}), $a=\frac{M+m}{2}y,$ $r=\frac{1}{2}\left( M-m\right)
\left\Vert y\right\Vert $ we get%
\begin{equation*}
\left\Vert x\right\Vert \left\Vert y\right\Vert \sqrt{Mm}\leq \frac{M+m}{2}%
\func{Re}\left\langle x,y\right\rangle
\end{equation*}%
giving 
\begin{equation*}
0\leq \left\Vert x\right\Vert \left\Vert y\right\Vert -\func{Re}\left\langle
x,y\right\rangle \leq \frac{\left( \sqrt{M}-\sqrt{m}\right) ^{2}}{2\sqrt{mM}}%
\func{Re}\left\langle x,y\right\rangle .
\end{equation*}%
Following the same arguments as in the proof of Proposition \ref{p2.4}, we
deduce the desired inequality (\ref{2.24}).
\end{proof}

For some results related to triangle inequality in inner product spaces, see 
\cite{JBDFTM}, \cite{SMK}, \cite{PMM} and \cite{DKR}.

\section{Some Gr\"{u}ss Type Inequalities\label{s4}}

We may state the following result.

\begin{theorem}
\label{t4.1}Let $\left( H;\left\langle \cdot ,\cdot \right\rangle \right) $
be an inner product space over the real or complex number field $\mathbb{K}$ 
$\left( \mathbb{K}=\mathbb{R},\mathbb{K}=\mathbb{C}\right) $ and $x,y,e\in H$
with $\left\Vert e\right\Vert =1.$ If $r_{1},r_{2}\in \left( 0,1\right) $
and 
\begin{equation}
\left\Vert x-e\right\Vert \leq r_{1},\ \ \ \ \left\Vert y-e\right\Vert \leq
r_{2},  \label{4.1}
\end{equation}%
then we have the inequality%
\begin{equation}
\left\vert \left\langle x,y\right\rangle -\left\langle x,e\right\rangle
\left\langle e,y\right\rangle \right\vert \leq r_{1}r_{2}\left\Vert
x\right\Vert \left\Vert y\right\Vert .  \label{4.2}
\end{equation}%
The inequality (\ref{4.2}) is sharp in the sense that the constant $c=1$ in
front of $r_{1}r_{2}$ cannot be replaced by a smaller constant.
\end{theorem}

\begin{proof}
Apply Schwarz's inequality in $\left( H;\left\langle \cdot ,\cdot
\right\rangle \right) $ for the vectors $x-\left\langle x,e\right\rangle e,$ 
$y-\left\langle y,e\right\rangle e,$ to get (see also \cite{SSD3})%
\begin{equation}
\left\vert \left\langle x,y\right\rangle -\left\langle x,e\right\rangle
\left\langle e,y\right\rangle \right\vert ^{2}\leq \left( \left\Vert
x\right\Vert ^{2}-\left\vert \left\langle x,e\right\rangle \right\vert
^{2}\right) \left( \left\Vert y\right\Vert ^{2}-\left\vert \left\langle
y,e\right\rangle \right\vert ^{2}\right) .  \label{4.3}
\end{equation}%
Using Theorem \ref{t2.1} for $a=e,$ we may state that%
\begin{equation}
\left\Vert x\right\Vert ^{2}-\left\vert \left\langle x,e\right\rangle
\right\vert ^{2}\leq r_{1}^{2}\left\Vert x\right\Vert ^{2},\ \ \ \ \ \
\left\Vert y\right\Vert ^{2}-\left\vert \left\langle y,e\right\rangle
\right\vert ^{2}\leq r_{2}^{2}\left\Vert y\right\Vert ^{2}.  \label{4.4}
\end{equation}%
Utilizing (\ref{4.3}) and (\ref{4.4}), we deduce%
\begin{equation}
\left\vert \left\langle x,y\right\rangle -\left\langle x,e\right\rangle
\left\langle e,y\right\rangle \right\vert ^{2}\leq
r_{1}^{2}r_{2}^{2}\left\Vert x\right\Vert ^{2}\left\Vert y\right\Vert ^{2},
\label{4.5}
\end{equation}%
which is clearly equivalent to the desired inequality (\ref{4.2}).

The sharpness of the constant follows by the fact that for $x=y,$ $%
r_{1}=r_{2}=r,$ we get from (\ref{4.2})%
\begin{equation}
\left\Vert x\right\Vert ^{2}-\left\vert \left\langle x,e\right\rangle
\right\vert ^{2}\leq r^{2}\left\Vert x\right\Vert ^{2}  \label{4.6}
\end{equation}%
provided $\left\Vert e\right\Vert =1$ and $\left\Vert x-e\right\Vert \leq
r<1.$ The inequality (\ref{4.6}) is sharp, as shown in Theorem \ref{t2.1},
and the theorem is thus proved.
\end{proof}

Another companion of the Gr\"{u}ss inequality may be stated as well.

\begin{theorem}
\label{t4.2}Let $\left( H;\left\langle \cdot ,\cdot \right\rangle \right) $
be an inner product space over $\mathbb{K}$ and $x,y,e\in H$ with $%
\left\Vert e\right\Vert =1.$ Suppose also that $a,A,b,B\in \mathbb{K}$ $%
\left( \mathbb{K}=\mathbb{R},\mathbb{C}\right) $ such that $\func{Re}\left( A%
\overline{a}\right) ,$ $\func{Re}\left( B\overline{b}\right) >0.$ If either%
\begin{equation}
\func{Re}\left\langle Ae-x,x-ae\right\rangle \geq 0,\ \ \func{Re}%
\left\langle Be-y,y-be\right\rangle \geq 0,\   \label{4.7}
\end{equation}%
or, equivalently,%
\begin{equation}
\left\Vert x-\frac{a+A}{2}e\right\Vert \leq \frac{1}{2}\left\vert
A-a\right\vert ,\ \ \left\Vert y-\frac{b+B}{2}e\right\Vert \leq \frac{1}{2}%
\left\vert B-b\right\vert ,  \label{4.8}
\end{equation}%
holds, then we have the inequality%
\begin{equation}
\left\vert \left\langle x,y\right\rangle -\left\langle x,e\right\rangle
\left\langle e,y\right\rangle \right\vert \leq \frac{1}{4}\cdot \frac{%
\left\vert A-a\right\vert \left\vert B-b\right\vert }{\sqrt{\func{Re}\left( A%
\overline{a}\right) \func{Re}\left( B\overline{b}\right) }}\left\vert
\left\langle x,e\right\rangle \left\langle e,y\right\rangle \right\vert .
\label{4.9}
\end{equation}%
The constant $\frac{1}{4}$ is best possible.
\end{theorem}

\begin{proof}
We know, by (\ref{4.3}), that%
\begin{equation}
\left\vert \left\langle x,y\right\rangle -\left\langle x,e\right\rangle
\left\langle e,y\right\rangle \right\vert ^{2}\leq \left( \left\Vert
x\right\Vert ^{2}-\left\vert \left\langle x,e\right\rangle \right\vert
^{2}\right) \left( \left\Vert y\right\Vert ^{2}-\left\vert \left\langle
y,e\right\rangle \right\vert ^{2}\right) .  \label{4.10}
\end{equation}%
If we use Corollary \ref{c2.3}, then we may state that%
\begin{equation}
\left\Vert x\right\Vert ^{2}-\left\vert \left\langle x,e\right\rangle
\right\vert ^{2}\leq \frac{1}{4}\cdot \frac{\left\vert A-a\right\vert ^{2}}{%
\func{Re}\left( A\overline{a}\right) }\left\vert \left\langle
x,e\right\rangle \right\vert ^{2}  \label{4.11}
\end{equation}%
and%
\begin{equation}
\left\Vert y\right\Vert ^{2}-\left\vert \left\langle y,e\right\rangle
\right\vert ^{2}\leq \frac{1}{4}\cdot \frac{\left\vert B-b\right\vert ^{2}}{%
\func{Re}\left( B\overline{b}\right) }\left\vert \left\langle
y,e\right\rangle \right\vert ^{2}.  \label{4.12}
\end{equation}%
Utilizing (\ref{4.10}) -- (\ref{4.12}), we deduce%
\begin{equation*}
\left\vert \left\langle x,y\right\rangle -\left\langle x,e\right\rangle
\left\langle e,y\right\rangle \right\vert ^{2}\leq \frac{1}{16}\cdot \frac{%
\left\vert A-a\right\vert ^{2}\left\vert B-b\right\vert ^{2}}{\func{Re}%
\left( A\overline{a}\right) \func{Re}\left( B\overline{b}\right) }\left\vert
\left\langle x,e\right\rangle \left\langle e,y\right\rangle \right\vert ^{2},
\end{equation*}%
which is clearly equivalent to the desired inequality (\ref{4.9}).

The sharpness of the constant follows from Corollary \ref{c2.3}, and we omit
the details.
\end{proof}

\begin{remark}
With the assumptions of Theorem \ref{t4.2} and if $\left\langle
x,e\right\rangle ,\left\langle y,e\right\rangle \neq 0$ (that is actually
the interesting case), one has the inequality%
\begin{equation}
\left\vert \frac{\left\langle x,y\right\rangle }{\left\langle
x,e\right\rangle \left\langle e,y\right\rangle }-1\right\vert \leq \frac{1}{4%
}\cdot \frac{\left\vert A-a\right\vert \left\vert B-b\right\vert }{\sqrt{%
\func{Re}\left( A\overline{a}\right) \func{Re}\left( B\overline{b}\right) }}.
\label{4.13}
\end{equation}%
The constant $\frac{1}{4}$ is best possible.
\end{remark}

\begin{remark}
The inequality (\ref{4.9}) provides a better bound for the quantity 
\begin{equation*}
\left\vert \left\langle x,y\right\rangle -\left\langle x,e\right\rangle
\left\langle e,y\right\rangle \right\vert
\end{equation*}%
than (2.3) of \cite{SSD3}.
\end{remark}

For some recent results on Gr\"{u}ss type inequalities in inner product
spaces, see \cite{SSD0}, \cite{SSD00} and \cite{PFR}.

\section{Reverses of Bessel's Inequality\label{s5}}

Let $\left( H;\left\langle \cdot ,\cdot \right\rangle \right) $ be a real or
complex infinite dimensional Hilbert space and $\left( e_{i}\right) _{i\in 
\mathbb{N}}$ an orthornormal family in $H$, i.e., we recall that $%
\left\langle e_{i},e_{j}\right\rangle =0$ if $i,j\in \mathbb{N}$, $i\neq j$
and $\left\Vert e_{i}\right\Vert =1$ for $i\in \mathbb{N}$.

It is well known that, if $x\in H,$ then the sum $\sum_{i=1}^{\infty
}\left\vert \left\langle x,e_{i}\right\rangle \right\vert ^{2}$ is
convergent and the following inequality, called \textit{Bessel's inequality}%
\begin{equation}
\sum_{i=1}^{\infty }\left\vert \left\langle x,e_{i}\right\rangle \right\vert
^{2}\leq \left\Vert x\right\Vert ^{2},  \label{5.1}
\end{equation}%
holds.

If $\ell ^{2}\left( \mathbb{K}\right) :=\left\{ \mathbf{a}=\left(
a_{i}\right) _{i\in \mathbb{N}}\subset \mathbb{K}\left\vert
\sum_{i=1}^{\infty }\left\vert a_{i}\right\vert ^{2}\right. <\infty \right\}
,$ where $\mathbb{K}=\mathbb{C}$ or $\mathbb{K}=\mathbb{R}$, is the Hilbert
space of all complex or real sequences that are $2$-summable and $\mathbf{%
\lambda }=\left( \lambda _{i}\right) _{i\in \mathbb{N}}\in \ell ^{2}\left( 
\mathbb{K}\right) ,$ then the sum $\sum_{i=1}^{\infty }\lambda _{i}e_{i}$ is
convergent in $H$ and if $y:=\sum_{i=1}^{\infty }\lambda _{i}e_{i}\in H,$
then $\left\Vert y\right\Vert =\left( \sum_{i=1}^{\infty }\left\vert \lambda
_{i}\right\vert ^{2}\right) ^{\frac{1}{2}}.$

We may state the following result.

\begin{theorem}
\label{t5.1}Let $\left( H;\left\langle \cdot ,\cdot \right\rangle \right) $
be an infinite dimensional Hilbert space over the real or complex number
field $\mathbb{K}$, $\left( e_{i}\right) _{i\in \mathbb{N}}$ an orthornormal
family in $H,$ $\mathbf{\lambda }=\left( \lambda _{i}\right) _{i\in \mathbb{N%
}}\in \ell ^{2}\left( \mathbb{K}\right) $ and $r>0$ with the property that%
\begin{equation}
\sum_{i=1}^{\infty }\left\vert \lambda _{i}\right\vert ^{2}>r^{2}.
\label{5.2}
\end{equation}%
If $x\in H$ is such that%
\begin{equation}
\left\Vert x-\sum_{i=1}^{\infty }\lambda _{i}e_{i}\right\Vert \leq r,
\label{5.3}
\end{equation}%
then we have the inequality%
\begin{align}
\left\Vert x\right\Vert ^{2}& \leq \frac{\left( \sum_{i=1}^{\infty }\func{Re}%
\left[ \overline{\lambda _{i}}\left\langle x,e_{i}\right\rangle \right]
\right) ^{2}}{\sum_{i=1}^{\infty }\left\vert \lambda _{i}\right\vert
^{2}-r^{2}}  \label{5.4} \\
& \leq \frac{\left\vert \sum_{i=1}^{\infty }\overline{\lambda _{i}}%
\left\langle x,e_{i}\right\rangle \right\vert ^{2}}{\sum_{i=1}^{\infty
}\left\vert \lambda _{i}\right\vert ^{2}-r^{2}}  \notag \\
& \leq \frac{\sum_{i=1}^{\infty }\left\vert \lambda _{i}\right\vert ^{2}}{%
\sum_{i=1}^{\infty }\left\vert \lambda _{i}\right\vert ^{2}-r^{2}}%
\sum_{i=1}^{\infty }\left\vert \left\langle x,e_{i}\right\rangle \right\vert
^{2};  \notag
\end{align}%
and%
\begin{eqnarray}
0 &\leq &\left\Vert x\right\Vert ^{2}-\sum_{i=1}^{\infty }\left\vert
\left\langle x,e_{i}\right\rangle \right\vert ^{2}  \label{5.5} \\
&\leq &\frac{r^{2}}{\sum_{i=1}^{\infty }\left\vert \lambda _{i}\right\vert
^{2}-r^{2}}\sum_{i=1}^{\infty }\left\vert \left\langle x,e_{i}\right\rangle
\right\vert ^{2}.
\end{eqnarray}
\end{theorem}

\begin{proof}
Applying the third inequality in (\ref{2.2}) for $a=\sum_{i=1}^{\infty
}\lambda _{i}e_{i}\in H,$ we have%
\begin{equation}
\left\Vert x\right\Vert ^{2}\left\Vert \sum_{i=1}^{\infty }\lambda
_{i}e_{i}\right\Vert ^{2}-\left[ \func{Re}\left\langle x,\sum_{i=1}^{\infty
}\lambda _{i}e_{i}\right\rangle \right] ^{2}\leq r^{2}\left\Vert
x\right\Vert ^{2}  \label{5.6}
\end{equation}%
and since%
\begin{align*}
\left\Vert \sum_{i=1}^{\infty }\lambda _{i}e_{i}\right\Vert ^{2}&
=\sum_{i=1}^{\infty }\left\vert \lambda _{i}\right\vert ^{2}, \\
\func{Re}\left\langle x,\sum_{i=1}^{\infty }\lambda _{i}e_{i}\right\rangle &
=\sum_{i=1}^{\infty }\func{Re}\left[ \overline{\lambda _{i}}\left\langle
x,e_{i}\right\rangle \right] ,
\end{align*}%
then by (\ref{5.6}) we deduce%
\begin{equation*}
\left\Vert x\right\Vert ^{2}\sum_{i=1}^{\infty }\left\vert \lambda
_{i}\right\vert ^{2}-\left[ \func{Re}\left\langle x,\sum_{i=1}^{\infty
}\lambda _{i}e_{i}\right\rangle \right] ^{2}\leq r^{2}\left\Vert
x\right\Vert ^{2},
\end{equation*}%
giving the first inequality in (\ref{5.4}).

The second inequality is obvious by the modulus property.

The last inequality follows by the Cauchy-Bunyakovsky-Schwarz inequality%
\begin{equation*}
\left\vert \sum_{i=1}^{\infty }\overline{\lambda _{i}}\left\langle
x,e_{i}\right\rangle \right\vert ^{2}\leq \sum_{i=1}^{\infty }\left\vert
\lambda _{i}\right\vert ^{2}\sum_{i=1}^{\infty }\left\vert \left\langle
x,e_{i}\right\rangle \right\vert ^{2}.
\end{equation*}%
The inequality (\ref{5.5}) follows by the last inequality in (\ref{5.4}) on
subtracting in both sides the quantity $\sum_{i=1}^{\infty }\left\vert
\left\langle x,e_{i}\right\rangle \right\vert ^{2}<\infty .$
\end{proof}

The following result provides a generalization for the reverse of Bessel's
inequality obtained in \cite{SSD6}.

\begin{theorem}
\label{t5.2}Let $\left( H;\left\langle \cdot ,\cdot \right\rangle \right) $
and $\left( e_{i}\right) _{i\in \mathbb{N}}$ be as in Theorem \ref{t5.1}.
Suppose that $\mathbf{\Gamma }=\left( \Gamma _{i}\right) _{i\in \mathbb{N}%
}\in \ell ^{2}\left( \mathbb{K}\right) ,$ $\mathbf{\gamma }=\left( \gamma
_{i}\right) _{i\in \mathbb{N}}\in \ell ^{2}\left( \mathbb{K}\right) $ are
sequences of real or complex numbers such that%
\begin{equation}
\sum_{i=1}^{\infty }\func{Re}\left( \Gamma _{i}\overline{\gamma _{i}}\right)
>0.  \label{5.7}
\end{equation}%
If $x\in H$ is such that either%
\begin{equation}
\left\Vert x-\sum_{i=1}^{\infty }\frac{\Gamma _{i}+\gamma _{i}}{2}%
e_{i}\right\Vert \leq \frac{1}{2}\left( \sum_{i=1}^{\infty }\left\vert
\Gamma _{i}-\gamma _{i}\right\vert ^{2}\right) ^{\frac{1}{2}}  \label{5.8}
\end{equation}%
or, equivalently,%
\begin{equation}
\func{Re}\left\langle \sum_{i=1}^{\infty }\Gamma
_{i}e_{i}-x,x-\sum_{i=1}^{\infty }\gamma _{i}e_{i}\right\rangle \geq 0
\label{5.9}
\end{equation}%
holds, then we have the inequalities%
\begin{align}
\left\Vert x\right\Vert ^{2}& \leq \frac{1}{4}\cdot \frac{\left(
\sum_{i=1}^{\infty }\func{Re}\left[ \left( \overline{\Gamma _{i}}+\overline{%
\gamma _{i}}\right) \left\langle x,e_{i}\right\rangle \right] \right) ^{2}}{%
\sum_{i=1}^{\infty }\func{Re}\left( \Gamma _{i}\overline{\gamma _{i}}\right) 
}  \label{5.10} \\
& \leq \frac{1}{4}\cdot \frac{\left\vert \sum_{i=1}^{\infty }\left( 
\overline{\Gamma _{i}}+\overline{\gamma _{i}}\right) \left\langle
x,e_{i}\right\rangle \right\vert ^{2}}{\sum_{i=1}^{\infty }\func{Re}\left(
\Gamma _{i}\overline{\gamma _{i}}\right) }  \notag \\
& \leq \frac{1}{4}\cdot \frac{\sum_{i=1}^{\infty }\left\vert \Gamma
_{i}+\gamma _{i}\right\vert ^{2}}{\sum_{i=1}^{\infty }\func{Re}\left( \Gamma
_{i}\overline{\gamma _{i}}\right) }\sum_{i=1}^{\infty }\left\vert
\left\langle x,e_{i}\right\rangle \right\vert ^{2}.  \notag
\end{align}%
The constant $\frac{1}{4}$ is best possible in all inequalities in (\ref%
{5.10}).

We also have the inequalities:%
\begin{equation}
0\leq \left\Vert x\right\Vert ^{2}-\sum_{i=1}^{\infty }\left\vert
\left\langle x,e_{i}\right\rangle \right\vert ^{2}\leq \frac{1}{4}\cdot 
\frac{\sum_{i=1}^{\infty }\left\vert \Gamma _{i}-\gamma _{i}\right\vert ^{2}%
}{\sum_{i=1}^{\infty }\func{Re}\left( \Gamma _{i}\overline{\gamma _{i}}%
\right) }\sum_{i=1}^{\infty }\left\vert \left\langle x,e_{i}\right\rangle
\right\vert ^{2}.  \label{5.11}
\end{equation}%
Here the constant $\frac{1}{4}$ is also best possible.
\end{theorem}

\begin{proof}
Since $\mathbf{\Gamma }$, $\mathbf{\gamma }\in \ell ^{2}\left( \mathbb{K}%
\right) ,$ then also $\frac{1}{2}\left( \mathbf{\Gamma }+\mathbf{\gamma }%
\right) \in \ell ^{2}\left( \mathbb{K}\right) ,$ showing that the series%
\begin{equation*}
\sum_{i=1}^{\infty }\left\vert \frac{\Gamma _{i}+\gamma _{i}}{2}\right\vert
^{2},\ \sum_{i=1}^{\infty }\left\vert \frac{\Gamma _{i}-\gamma _{i}}{2}%
\right\vert ^{2}\text{ and}\ \sum_{i=1}^{\infty }\func{Re}\left( \Gamma _{i}%
\overline{\gamma _{i}}\right) 
\end{equation*}%
are convergent. Also, the series 
\begin{equation*}
\sum_{i=1}^{\infty }\Gamma _{i}e_{i},\text{ }\sum_{i=1}^{\infty }\gamma
_{i}e_{i}\text{ and }\sum_{i=1}^{\infty }\frac{\gamma _{i}+\Gamma _{i}}{2}%
e_{i}
\end{equation*}
are convergent in the Hilbert space $H.$

The equivalence of the conditions (\ref{5.8}) and (\ref{5.9}) follows by the
fact that in an inner product space we have, for $x,z,Z\in H,$ $\func{Re}%
\left\langle Z-x,x-z\right\rangle \geq 0$ is equivalent to $\left\Vert x-%
\frac{z+Z}{2}\right\Vert \leq \frac{1}{2}\left\Vert Z-z\right\Vert ,$ and we
omit the details.

Now, we observe that the inequalities (\ref{5.10}) and (\ref{5.11}) follow
from Theorem \ref{t5.1} on choosing $\lambda _{i}=\frac{\gamma _{i}+\Gamma
_{i}}{2},$ $i\in \mathbb{N}$ and $r=\frac{1}{2}\left( \sum_{i=1}^{\infty
}\left\vert \Gamma _{i}-\gamma _{i}\right\vert ^{2}\right) ^{\frac{1}{2}}.$

The fact that $\frac{1}{4}$ is the best constant in both (\ref{5.10}) and (%
\ref{5.11}) follows from Theorem \ref{t2.2} and Corollary \ref{c2.3}, and we
omit the details.
\end{proof}

\begin{remark}
Note that (\ref{5.10}) improves (\ref{1.17}) and (\ref{5.11}) improves (\ref%
{1.18}), that have been obtained in \cite{SSD6}.
\end{remark}

For some recent results related to Bessel inequality, see \cite{SSD01}, \cite%
{SSDJS}, \cite{HXC} and  \cite{GH1}.

\section{Some Gr\"{u}ss Type Inequalities for Orthonormal Families\label{s6}}

The following result related to Gr\"{u}ss inequality in inner product
spaces, holds.

\begin{theorem}
\label{t6.1}Let $\left( H;\left\langle \cdot ,\cdot \right\rangle \right) $
be an infinite dimensional Hilbert space over the real or complex number
field $\mathbb{K}$, and $\left( e_{i}\right) _{i\in \mathbb{N}}$ an
orthornormal family in $H.$ Assume that $\mathbf{\lambda }=\left( \lambda
_{i}\right) _{i\in \mathbb{N}},\ \mathbf{\mu }=\left( \mu _{i}\right) _{i\in 
\mathbb{N}}\in \ell ^{2}\left( \mathbb{K}\right) $ and $r_{1},r_{2}>0$ with
the properties that%
\begin{equation}
\sum_{i=1}^{\infty }\left\vert \lambda _{i}\right\vert ^{2}>r_{1}^{2},\ \ \
\sum_{i=1}^{\infty }\left\vert \mu _{i}\right\vert ^{2}>r_{2}^{2}.
\label{6.1}
\end{equation}%
If $x,y\in H$ are such that%
\begin{equation}
\left\Vert x-\sum_{i=1}^{\infty }\lambda _{i}e_{i}\right\Vert \leq r_{1},\ \
\ \ \ \ \left\Vert y-\sum_{i=1}^{\infty }\mu _{i}e_{i}\right\Vert \leq r_{2},
\label{6.2}
\end{equation}%
then we have the inequalities%
\begin{eqnarray}
&&\left\vert \left\langle x,y\right\rangle -\sum_{i=1}^{\infty }\left\langle
x,e_{i}\right\rangle \left\langle e_{i},y\right\rangle \right\vert 
\label{6.3} \\
&\leq &\frac{r_{1}r_{2}}{\sqrt{\sum_{i=1}^{\infty }\left\vert \lambda
_{i}\right\vert ^{2}-r_{1}^{2}}\sqrt{\sum_{i=1}^{\infty }\left\vert \mu
_{i}\right\vert ^{2}-r_{2}^{2}}}\cdot \sqrt{\sum_{i=1}^{\infty }\left\vert
\left\langle x,e_{i}\right\rangle \right\vert ^{2}\sum_{i=1}^{\infty
}\left\vert \left\langle y,e_{i}\right\rangle \right\vert ^{2}}  \notag \\
&\leq &\frac{r_{1}r_{2}\left\Vert x\right\Vert \left\Vert y\right\Vert }{%
\sqrt{\sum_{i=1}^{\infty }\left\vert \lambda _{i}\right\vert ^{2}-r_{1}^{2}}%
\sqrt{\sum_{i=1}^{\infty }\left\vert \mu _{i}\right\vert ^{2}-r_{2}^{2}}}. 
\notag
\end{eqnarray}
\end{theorem}

\begin{proof}
Applying Schwarz's inequality for the vectors $x-\sum_{i=1}^{\infty
}\left\langle x,e_{i}\right\rangle e_{i},$ $y-\sum_{i=1}^{\infty
}\left\langle y,e_{i}\right\rangle e_{i},$ we have%
\begin{multline}
\left\vert \left\langle x-\sum_{i=1}^{\infty }\left\langle
x,e_{i}\right\rangle e_{i},y-\sum_{i=1}^{\infty }\left\langle
y,e_{i}\right\rangle e_{i}\right\rangle \right\vert ^{2}  \label{6.4} \\
\leq \left\Vert x-\sum_{i=1}^{\infty }\left\langle x,e_{i}\right\rangle
e_{i}\right\Vert ^{2}\left\Vert y-\sum_{i=1}^{\infty }\left\langle
y,e_{i}\right\rangle e_{i}\right\Vert ^{2}.
\end{multline}%
Since%
\begin{equation*}
\left\langle x-\sum_{i=1}^{\infty }\left\langle x,e_{i}\right\rangle
e_{i},y-\sum_{i=1}^{\infty }\left\langle y,e_{i}\right\rangle
e_{i}\right\rangle =\left\langle x,y\right\rangle -\sum_{i=1}^{\infty
}\left\langle x,e_{i}\right\rangle \left\langle e_{i},y\right\rangle 
\end{equation*}%
and%
\begin{equation*}
\left\Vert x-\sum_{i=1}^{\infty }\left\langle x,e_{i}\right\rangle
e_{i}\right\Vert ^{2}=\left\Vert x\right\Vert ^{2}-\sum_{i=1}^{\infty
}\left\vert \left\langle x,e_{i}\right\rangle \right\vert ^{2},
\end{equation*}%
then by (\ref{5.5}) applied for $x$ and $y,$ and from (\ref{6.4}), we deduce
the first part of (\ref{6.3}).

The second part follows by Bessel's inequality.
\end{proof}

The following Gr\"{u}ss type inequality may be stated as well.

\begin{theorem}
\label{t6.2}Let $\left( H;\left\langle \cdot ,\cdot \right\rangle \right) $
be an infinite dimensional Hilbert space and $\left( e_{i}\right) _{i\in 
\mathbb{N}}$ an orthornormal family in $H.$ Suppose that $\left( \Gamma
_{i}\right) _{i\in \mathbb{N}},$ $\left( \gamma _{i}\right) _{i\in \mathbb{N}%
},$ $\left( \phi _{i}\right) _{i\in \mathbb{N}},$ $\left( \Phi _{i}\right)
_{i\in \mathbb{N}}\in \ell ^{2}\left( \mathbb{K}\right) $ are sequences of
real and complex numbers such that%
\begin{equation}
\sum_{i=1}^{\infty }\func{Re}\left( \Gamma _{i}\overline{\gamma _{i}}\right)
>0,\ \ \ \sum_{i=1}^{\infty }\func{Re}\left( \Phi _{i}\overline{\phi _{i}}%
\right) >0.  \label{6.5}
\end{equation}%
If $x,y\in H$ are such that either%
\begin{align}
\left\Vert x-\sum_{i=1}^{\infty }\frac{\Gamma _{i}+\gamma _{i}}{2}\cdot
e_{i}\right\Vert & \leq \frac{1}{2}\left( \sum_{i=1}^{\infty }\left\vert
\Gamma _{i}-\gamma _{i}\right\vert ^{2}\right) ^{\frac{1}{2}}  \label{6.6} \\
\left\Vert y-\sum_{i=1}^{\infty }\frac{\Phi _{i}+\phi _{i}}{2}\cdot
e_{i}\right\Vert & \leq \frac{1}{2}\left( \sum_{i=1}^{\infty }\left\vert
\Phi _{i}-\phi _{i}\right\vert ^{2}\right) ^{\frac{1}{2}}  \notag
\end{align}%
or, equivalently,%
\begin{align}
\func{Re}\left\langle \sum_{i=1}^{\infty }\Gamma
_{i}e_{i}-x,x-\sum_{i=1}^{\infty }\gamma _{i}e_{i}\right\rangle & \geq 0,
\label{6.7} \\
\func{Re}\left\langle \sum_{i=1}^{\infty }\Phi
_{i}e_{i}-y,y-\sum_{i=1}^{\infty }\phi _{i}e_{i}\right\rangle & \geq 0, 
\notag
\end{align}%
holds, then we have the inequality%
\begin{align}
& \left\vert \left\langle x,y\right\rangle -\sum_{i=1}^{\infty }\left\langle
x,e_{i}\right\rangle \left\langle e_{i},y\right\rangle \right\vert 
\label{6.8} \\
& \leq \frac{1}{4}\cdot \frac{\left( \sum_{i=1}^{\infty }\left\vert \Gamma
_{i}-\gamma _{i}\right\vert ^{2}\right) ^{\frac{1}{2}}\left(
\sum_{i=1}^{\infty }\left\vert \Phi _{i}-\phi _{i}\right\vert ^{2}\right) ^{%
\frac{1}{2}}}{\left( \sum_{i=1}^{\infty }\func{Re}\left( \Gamma _{i}%
\overline{\gamma _{i}}\right) \right) ^{\frac{1}{2}}\left(
\sum_{i=1}^{\infty }\func{Re}\left( \Phi _{i}\overline{\phi _{i}}\right)
\right) ^{\frac{1}{2}}}  \notag \\
& \times \left( \sum_{i=1}^{\infty }\left\vert \left\langle
x,e_{i}\right\rangle \right\vert ^{2}\right) ^{\frac{1}{2}}\left(
\sum_{i=1}^{\infty }\left\vert \left\langle y,e_{i}\right\rangle \right\vert
^{2}\right) ^{\frac{1}{2}}  \notag \\
& \leq \frac{1}{4}\cdot \frac{\left( \sum_{i=1}^{\infty }\left\vert \Gamma
_{i}-\gamma _{i}\right\vert ^{2}\right) ^{\frac{1}{2}}\left(
\sum_{i=1}^{\infty }\left\vert \Phi _{i}-\phi _{i}\right\vert ^{2}\right) ^{%
\frac{1}{2}}}{\left[ \sum_{i=1}^{\infty }\func{Re}\left( \Gamma _{i}%
\overline{\gamma _{i}}\right) \right] ^{\frac{1}{2}}\left[
\sum_{i=1}^{\infty }\func{Re}\left( \Phi _{i}\overline{\phi _{i}}\right) %
\right] ^{\frac{1}{2}}}\left\Vert x\right\Vert \left\Vert y\right\Vert  
\notag
\end{align}%
The constant $\frac{1}{4}$ is best possible in the first inequality.
\end{theorem}

\begin{proof}
Follows by (\ref{5.11}) and (\ref{6.4}).

The best constant follows from Theorem \ref{t4.2}, and we omit the details.
\end{proof}

\begin{remark}
We note that the inequality (\ref{6.8}) is better than the inequality (3.3)
in \cite{SSD6}. We omit the details.
\end{remark}

\section{Integral Inequalities\label{s7}}

Let $\left( \Omega ,\Sigma ,\mu \right) $ be a measurable space consisting
of a set $\Omega ,$ a $\sigma -$algebra of parts $\Sigma $ and a countably
additive and positive measure $\mu $ on $\Sigma $ with values $\mathbb{R\cup 
}\left\{ \infty \right\} .$ Let $\rho \geq 0$ be a $g-$measurable function
on $\Omega $ with $\int_{\Omega }\rho \left( s\right) d\mu \left( s\right)
=1.$ Denote by $L_{\rho }^{2}\left( \Omega ,\mathbb{K}\right) $ the Hilbert
space of all real or complex valued functions defined on $\Omega $ and $%
2-\rho -$integrable on $\Omega ,$ i.e.,%
\begin{equation}
\int_{\Omega }\rho \left( s\right) \left\vert f\left( s\right) \right\vert
^{2}d\mu \left( s\right) <\infty .  \label{7.1}
\end{equation}%
It is obvious that the following inner product%
\begin{equation}
\left\langle f,g\right\rangle _{\rho }:=\int_{\Omega }\rho \left( s\right)
f\left( s\right) \overline{g\left( s\right) }d\mu \left( s\right) ,
\label{7.2}
\end{equation}%
generates the norm $\left\Vert f\right\Vert _{\rho }:=\left( \int_{\Omega
}\rho \left( s\right) \left\vert f\left( s\right) \right\vert ^{2}d\mu
\left( s\right) \right) ^{\frac{1}{2}}$ of $L_{\rho }^{2}\left( \Omega ,%
\mathbb{K}\right) ,$ and all the above results may be stated for integrals.

It is important to observe that, if 
\begin{equation}
\func{Re}\left[ f\left( s\right) \overline{g\left( s\right) }\right] \geq 0%
\text{ \ for }\mu -\text{a.e. }s\in \Omega ,  \label{7.3}
\end{equation}%
then, obviously,%
\begin{align}
\func{Re}\left\langle f,g\right\rangle _{\rho }& =\func{Re}\left[
\int_{\Omega }\rho \left( s\right) f\left( s\right) \overline{g\left(
s\right) }d\mu \left( s\right) \right]  \label{7.4} \\
& =\int_{\Omega }\rho \left( s\right) \func{Re}\left[ f\left( s\right) 
\overline{g\left( s\right) }\right] d\mu \left( s\right) \geq 0.  \notag
\end{align}%
The reverse is evidently not true in general.

Moreover, if the space is real, i.e., $\mathbb{K=R}$, then a sufficient
condition for (\ref{7.4}) to hold is:%
\begin{equation}
f\left( s\right) \geq 0,\ \ g\left( s\right) \geq 0\text{ \ for }\mu -\text{%
a.e. }s\in \Omega .  \label{7.5}
\end{equation}

We provide now, by the use of certain result obtained in Section \ref{s2},
some integral inequalities that may be used in practical applications.

\begin{proposition}
\label{p7.1}Let $f,g\in L_{\rho }^{2}\left( \Omega ,\mathbb{K}\right) $ and $%
r>0$ with the properties that%
\begin{equation}
\left\vert f\left( s\right) -g\left( s\right) \right\vert \leq r\leq
\left\vert g\left( s\right) \right\vert \ \text{\ for }\mu -\text{a.e. }s\in
\Omega ,  \label{7.6}
\end{equation}%
and $\int_{\Omega }\rho \left( s\right) \left\vert g\left( s\right)
\right\vert ^{2}d\mu \left( s\right) \neq r.$ Then we have the inequalities%
\begin{align}
0& \leq \int_{\Omega }\rho \left( s\right) \left\vert f\left( s\right)
\right\vert ^{2}d\mu \left( s\right) \int_{\Omega }\rho \left( s\right)
\left\vert g\left( s\right) \right\vert ^{2}d\mu \left( s\right) -\left\vert
\int_{\Omega }\rho \left( s\right) f\left( s\right) \overline{g\left(
s\right) }d\mu \left( s\right) \right\vert ^{2}  \label{7.7} \\
& \leq \int_{\Omega }\rho \left( s\right) \left\vert f\left( s\right)
\right\vert ^{2}d\mu \left( s\right) \int_{\Omega }\rho \left( s\right)
\left\vert g\left( s\right) \right\vert ^{2}d\mu \left( s\right)  \notag \\
& \ \ \ \ \ \ \ \ \ \ \ \ \ \ \ \ \ \ \ \ \ \ \ \ \ \ \ -\left[ \int_{\Omega
}\rho \left( s\right) \func{Re}\left( f\left( s\right) \overline{g\left(
s\right) }\right) d\mu \left( s\right) \right] ^{2}  \notag \\
& \leq r^{2}\int_{\Omega }\rho \left( s\right) \left\vert g\left( s\right)
\right\vert ^{2}d\mu \left( s\right) .  \notag
\end{align}%
The constant $c=1$ in front of $r^{2}$ is best possible.
\end{proposition}

The proof follows by Theorem \ref{t2.1} and we omit the details.

\begin{proposition}
\label{p7.2}Let $f,g\in L_{\rho }^{2}\left( \Omega ,\mathbb{K}\right) $ and $%
\gamma ,\Gamma \in \mathbb{K}$ such that $\func{Re}\left( \Gamma \overline{%
\gamma }\right) >0$ and%
\begin{equation}
\func{Re}\left[ \left( \Gamma g\left( s\right) -f\left( s\right) \right)
\left( \overline{f\left( s\right) }-\overline{\gamma }\overline{g\left(
s\right) }\right) \right] \geq 0\text{ \ for }\mu -\text{a.e. }s\in \Omega .
\label{7.8}
\end{equation}%
Then we have the inequalities%
\begin{align}
& \int_{\Omega }\rho \left( s\right) \left\vert f\left( s\right) \right\vert
^{2}d\mu \left( s\right) \int_{\Omega }\rho \left( s\right) \left\vert
g\left( s\right) \right\vert ^{2}d\mu \left( s\right)   \label{7.9} \\
& \leq \frac{1}{4}\cdot \frac{\left\{ \func{Re}\left[ \left( \overline{%
\Gamma }+\overline{\gamma }\right) \int_{\Omega }\rho \left( s\right)
f\left( s\right) \overline{g\left( s\right) }d\mu \left( s\right) \right]
\right\} ^{2}}{\func{Re}\left( \Gamma \overline{\gamma }\right) }  \notag \\
& \leq \frac{1}{4}\cdot \frac{\left\vert \Gamma +\gamma \right\vert ^{2}}{%
\func{Re}\left( \Gamma \overline{\gamma }\right) }\left\vert \int_{\Omega
}\rho \left( s\right) f\left( s\right) \overline{g\left( s\right) }d\mu
\left( s\right) \right\vert ^{2}.  \notag
\end{align}%
The constant $\frac{1}{4}$ is best possible in both inequalities.
\end{proposition}

The proof follows by Theorem \ref{t2.2} and we omit the details.

\begin{corollary}
\label{c7.3}With the assumptions of Proposition \ref{p7.2}, we have the
inequality%
\begin{align}
0& \leq \int_{\Omega }\rho \left( s\right) \left\vert f\left( s\right)
\right\vert ^{2}d\mu \left( s\right) \int_{\Omega }\rho \left( s\right)
\left\vert g\left( s\right) \right\vert ^{2}d\mu \left( s\right) 
\label{7.10} \\
& \ \ \ \ \ \ \ \ \ \ \ \ \ \ \ \ \ \ \ \ \ \ \ \ \ -\left\vert \int_{\Omega
}\rho \left( s\right) f\left( s\right) \overline{g\left( s\right) }d\mu
\left( s\right) \right\vert ^{2}  \notag \\
& \leq \frac{1}{4}\cdot \frac{\left\vert \Gamma -\gamma \right\vert ^{2}}{%
\func{Re}\left( \Gamma \overline{\gamma }\right) }\left\vert \int_{\Omega
}\rho \left( s\right) f\left( s\right) \overline{g\left( s\right) }d\mu
\left( s\right) \right\vert ^{2}.  \notag
\end{align}%
The constant $\frac{1}{4}$ is best possible.
\end{corollary}

\begin{remark}
If the space is real and we assume, for $M>m>0,$ that%
\begin{equation}
mg\left( s\right) \leq f\left( s\right) \leq Mg\left( s\right) \text{ \ for }%
\mu -\text{a.e. }s\in \Omega ,  \label{7.11}
\end{equation}%
then, by (\ref{7.9}) and (\ref{7.10}), we deduce the inequalities%
\begin{multline}
\int_{\Omega }\rho \left( s\right) \left[ f\left( s\right) \right] ^{2}d\mu
\left( s\right) \int_{\Omega }\rho \left( s\right) \left[ g\left( s\right) %
\right] ^{2}d\mu \left( s\right)   \label{7.12} \\
\leq \frac{1}{4}\cdot \frac{\left( M+m\right) ^{2}}{mM}\left[ \int_{\Omega
}\rho \left( s\right) f\left( s\right) g\left( s\right) d\mu \left( s\right) %
\right] ^{2}.
\end{multline}%
and 
\begin{align}
0& \leq \int_{\Omega }\rho \left( s\right) \left[ f\left( s\right) \right]
^{2}d\mu \left( s\right) \int_{\Omega }\rho \left( s\right) \left[ g\left(
s\right) \right] ^{2}d\mu \left( s\right)   \label{7.13} \\
& \ \ \ \ \ \ \ \ \ \ \ \ \ \ \ \ \ \ \ -\left[ \int_{\Omega }\rho \left(
s\right) f\left( s\right) g\left( s\right) d\mu \left( s\right) \right] ^{2}
\notag \\
& \leq \frac{1}{4}\cdot \frac{\left( M-m\right) ^{2}}{mM}\left[ \int_{\Omega
}\rho \left( s\right) f\left( s\right) g\left( s\right) d\mu \left( s\right) %
\right] ^{2}.  \notag
\end{align}%
The inequality (\ref{7.12}) is known in the literature as Cassel's
inequality.
\end{remark}

The following Gr\"{u}ss type integral inequality for real or complex-valued
functions also holds.

\begin{proposition}
\label{p.7.3}Let $f,g,h\in L_{\rho }^{2}\left( \Omega ,\mathbb{K}\right) $
with $\int_{\Omega }\rho \left( s\right) \left\vert h\left( s\right)
\right\vert ^{2}d\mu \left( s\right) =1$ and $a,A,b,B\in \mathbb{K}$ such
that $\func{Re}\left( A\overline{a}\right) ,\func{Re}\left( B\overline{b}%
\right) >0$ and%
\begin{eqnarray*}
\func{Re}\left[ \left( Ah\left( s\right) -f\left( s\right) \right) \left( 
\overline{f\left( s\right) }-\overline{a}\overline{h\left( s\right) }\right) %
\right]  &\geq &0, \\
\func{Re}\left[ \left( Ah\left( s\right) -g\left( s\right) \right) \left( 
\overline{g\left( s\right) }-\overline{b}\overline{h\left( s\right) }\right) %
\right]  &\geq &0\text{,}
\end{eqnarray*}%
for $\mu -$a.e. $s\in \Omega .$ Then we have the inequalities%
\begin{align}
& \left\vert \int_{\Omega }\rho \left( s\right) f\left( s\right) \overline{%
g\left( s\right) }d\mu \left( s\right) -\int_{\Omega }\rho \left( s\right)
f\left( s\right) \overline{h\left( s\right) }d\mu \left( s\right)
\int_{\Omega }\rho \left( s\right) h\left( s\right) \overline{g\left(
s\right) }d\mu \left( s\right) \right\vert  \\
& \leq \frac{1}{4}\cdot \frac{\left\vert A-a\right\vert \left\vert
B-b\right\vert }{\sqrt{\func{Re}\left( A\overline{a}\right) \func{Re}\left( B%
\overline{b}\right) }}\left\vert \int_{\Omega }\rho \left( s\right) f\left(
s\right) \overline{h\left( s\right) }d\mu \left( s\right) \int_{\Omega }\rho
\left( s\right) h\left( s\right) \overline{g\left( s\right) }d\mu \left(
s\right) \right\vert   \notag
\end{align}%
The constant $\frac{1}{4}$ is best possible.
\end{proposition}

The proof follows by Theorem \ref{t4.2}.

\begin{remark}
All the other inequalities in Sections \ref{s3} -- \ref{s6} may be used in a
similar way to obtain the corresponding integral inequalities. We omit the
details.
\end{remark}

\end{document}